\documentclass[pdflatex,sn-mathphys-num]{sn-jnl}

\usepackage{graphicx}%
\usepackage{multirow}%
\usepackage{amsmath,amssymb,amsfonts}%
\usepackage{amsthm}%
\usepackage{mathrsfs}%
\usepackage[title]{appendix}%
\usepackage{xcolor}%
\usepackage{textcomp}%
\usepackage{manyfoot}%
\usepackage{booktabs}%
\usepackage{algorithm}%
\usepackage{algorithmicx}%
\usepackage{algpseudocode}%
\usepackage{listings}%
\usepackage{subcaption}

\newcommand{\T}{\mathcal T}
\newcommand{\D}{\Delta_{\pi}^{\downarrow}}
\newcommand{\Fix}{\operatorname{Fix}}

\newcommand{\sortop}{\operatorname{sort}}

\theoremstyle{thmstyleone}%
\newtheorem{theorem}{Theorem}
\newtheorem{proposition}[theorem]{Proposition}%
\newtheorem{lemma}[theorem]{Lemma}
\newtheorem{corollary}[theorem]{Corollary}
\theoremstyle{definition}

\theoremstyle{thmstyletwo}%

\theoremstyle{thmstylethree}%

\begin{document}

\title[Largest-Angle Bisection and Statistical Regularization of Triangle Shapes]{Largest-Angle Bisection and Statistical Regularization of Triangle Shapes}


\author{\fnm{Jérôme} \sur{Michaud}}\email{jerome.michaud@mdu.se}



\affil{\orgdiv{Department of Business and Mathematics}, \orgname{Mälardalen University}, \orgaddress{\street{Universitetsplan 1}, \city{Västerås}, \postcode{722 20}, \country{Sweden}}}




\abstract{Largest-angle bisection (LAB) repeatedly divides a triangle by bisecting its largest angle. Although this procedure is known to preserve nondegenerate triangle shapes, less is known about the distribution of shapes across the resulting refinement tree. We show that its long-term behavior is governed by a distinction between pure and mixed branches. A branch that eventually selects only one of the two children may approach an exceptional family of shapes. By contrast, every branch that changes children infinitely often converges in shape to a right-isosceles triangle. This follows from a projective Lyapunov function that is non-increasing at every refinement and decreases strictly whenever the selected child changes. Consequently, a random branch converges to the right-isosceles shape with probability one, and the proportion of triangles near that shape tends to one across complete refinement levels. Largest-angle bisection therefore produces statistical regularization even though exceptional branches remain present.}




\maketitle

\section{Introduction}

Controlled mesh refinement is fundamental to the reliability of finite element computation. A refinement procedure should reduce element diameters while preserving element quality, typically expressed in two dimensions through a positive lower bound on triangle angles known as the minimum angle condition. Such geometric conditions underlie classical interpolation and finite element estimates; see, for example, \cite{BabuskaAziz1976,Ciarlet2002}.

Among triangular refinement rules, longest-edge bisection is especially well-studied. Rosenberg and Stenger established angle nondegeneracy \cite{RosenbergStenger1975}, while subsequent work proved mesh-size convergence and described the similarity classes generated by repeated refinement \cite{Stynes1980,adler1983bisection,Kearfott1978APO}. Perdomo and Plaza later formulated aspects of this behavior in a normalized shape space equipped with a hyperbolic metric \cite{PerdomoPlaza2014} and this space has been used to better understand the longterm behaviour of the iteration \cite{kalmanovich2026stability}. Longest-edge and related bisection procedures continue to play an important role in adaptive mesh generation; see also \cite{Rivara1984}.

Largest-angle bisection (LAB) instead selects a triangle's largest angle and joins its vertex to the opposite side along the angle bisector. Ismailescu, Kim, Kim, and Lee proved that repeated LAB produces shrinking mesh diameters and preserves a positive lower bound on the angles. They also showed that, except when the initial triangle is right-isosceles, the procedure generates infinitely many similarity classes \cite{ismailescu2019largestanglebisectionprocedure}. These results have been extended to largest-angle $n$-section in \cite{michaudKorotov2026}. Together, these results establish deterministic regularity, but they leave open a different question: how are triangle shapes distributed throughout the refinement tree? In particular, the existence of infinitely many similarity classes does not reveal whether these classes remain broadly distributed or whether most descendants become concentrated near a preferred shape.

Here we show that LAB exhibits a strong form of asymptotic shape selection. The two children produced at each bisection define two transformations of triangle shape space. Repetition of only one transformation leads to exceptional families of limiting shapes. Branches that continue to change between the two transformations behave differently: their shapes converge to the right-isosceles triangle. The mechanism is captured by a projective Lyapunov function that cannot increase under refinement and decreases strictly whenever the selected child changes. Thus, the exceptional behavior is associated precisely with branches that eventually repeat the same choice indefinitely.

This distinction between pure and mixed branches has both probabilistic and full-tree consequences. Under any independent random choice that assigns positive probability to both children, an infinite branch changes children infinitely often with probability one and therefore converges to the right-isosceles shape. Equivalently, in the complete refinement tree, the proportion of triangles lying outside any prescribed neighborhood of the right-isosceles shape tends to zero as the refinement level increases. LAB therefore displays \emph{statistical regularization}: exceptional branches remain present, but they become negligible among all descendants. This conclusion goes beyond deterministic angle bounds by describing the dominant geometry of the refinement tree rather than merely excluding degeneration.

The paper is organized as follows. Section~\ref{sec:model} introduces the projective shape space and the transformations induced by the two LAB children. Section~\ref{sec:dynamics} analyzes pure and mixed dynamics using a projective Lyapunov function and classifies the possible infinite branches. Section~\ref{sec:statistics} derives the probabilistic and complete-tree regularization results. Section~\ref{sec:conclusion} discusses the implications for mesh refinement and several directions for further research.

\section{Angle-projective shape space}\label{sec:model}

To study how triangle shapes are distributed under largest-angle bisection, we first need a normalized space in which similar triangles are represented by the same point. Perdomo and Plaza used such a shape-space viewpoint for longest-edge bisection (LEB), trisection and other $n$-sections, representing normalized triangles through side-length data and equipping the resulting space with a hyperbolic metric \cite{PerdomoPlaza2014,Perdomo2013,PerdomoPlaza2012}. Their construction has been used in \cite{kalmanovich2026stability} to study the distribution of the triangle shapes under LEB refinement and motivates the dynamical approach adopted here. For LAB, however, coordinates based directly on angles are more natural: the angle selected for refinement is determined by an angular comparison, and the two children can be described explicitly in terms of the parent angles.

We therefore introduce a projective angle space for triangle shapes. Ordering the angles removes the ambiguity caused by relabelling the vertices, while passing to ratios removes the common scale of a positive angle triple. In the resulting two-dimensional chart, the two LAB child operations become piecewise projective transformations. The piecewise character arises because the child angles must be reordered after each bisection.

\subsection{Angle-normalized shape space}

\begin{figure}
\centering
\begin{subfigure}[][][t]{0.48\textwidth}
    \includegraphics[width=\textwidth]{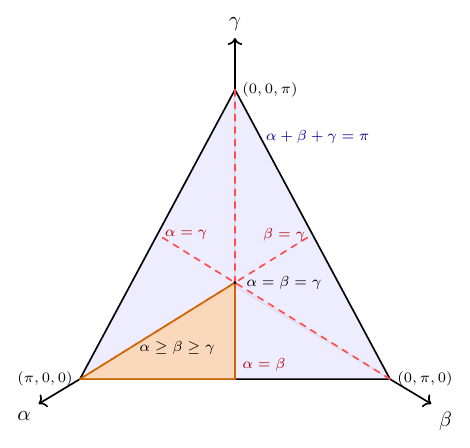}
    \caption{The angle simplex \(\alpha+\beta+\gamma=\pi\), viewed as a slice of the positive cone. The equality planes divide it into six sorting chambers, and the highlighted chamber represents \(\alpha\geq\beta\geq\gamma\).}
    \label{fig:first}
\end{subfigure}
\hfill
\begin{subfigure}[][][t]{0.48\textwidth}
    \includegraphics[width=\textwidth]{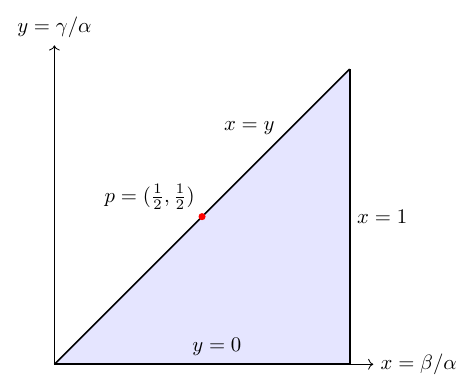}
    \caption{The projectively normalized shape domain \(\mathcal T=\{1\geq x\geq y>0\}\). The point \(p=(1/2,1/2)\) represents the right-isosceles triangle.}
    \label{fig:second}
\end{subfigure}

\caption{Illustration of the angle-projective shape space.}
\label{fig:angle-domain}
\end{figure}

Let a nondegenerate triangle have ordered angles
\[
\alpha\geq\beta\geq\gamma>0,
\qquad
\alpha+\beta+\gamma=\pi.
\]
Two triangles are similar if and only if their angle triples agree up to permutation. We remove this permutation ambiguity by arranging the angles in decreasing order and defining the ordered angle chamber
\[
\D=
\{(\alpha,\beta,\gamma):
\alpha\geq\beta\geq\gamma>0,
\ \alpha+\beta+\gamma=\pi\}.
\]
Thus \(\D\) contains one canonical representative of each unordered similarity class of nondegenerate triangles, see Figure~\ref{fig:first}.

The condition \(\alpha+\beta+\gamma=\pi\) already selects a representative from each positive ray of angle triples. Nevertheless, the action of LAB is more transparent when these triples are viewed projectively. We therefore divide the two smaller angles by the largest and introduce the chart
\begin{equation}\label{eq:chart}
\chi(\alpha,\beta,\gamma)
=
\left(\frac{\beta}{\alpha},\frac{\gamma}{\alpha}\right)
=(x,y).
\end{equation}
It maps \(\D\) bijectively onto
\[
\T=\{(x,y):1\geq x\geq y>0\},
\]
with inverse
\begin{equation}\label{eq:chartinverse}
\alpha=\frac{\pi}{1+x+y},
\qquad
\beta=\frac{\pi x}{1+x+y},
\qquad
\gamma=\frac{\pi y}{1+x+y}.
\end{equation}
The inequalities defining \(\T\) encode the ordering of the angles: the boundary \(x=1\) corresponds to \(\alpha=\beta\), the boundary \(x=y\) to \(\beta=\gamma\), and the excluded boundary \(y=0\) to degenerate angle triples. Hence \(\T\) provides one normalized point for every nondegenerate triangle shape, see Figure~\ref{fig:second}. The point
\[
p=\left(\frac12,\frac12\right)
\]
represents the right-isosceles triangle, whose ordered angles are \((\pi/2,\pi/4,\pi/4)\).

\subsection{LAB in angle coordinates}

\begin{figure}[ht]
\centering
\includegraphics[width=.7\textwidth]{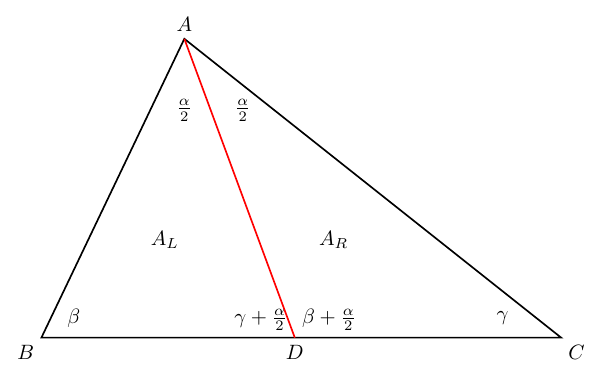}
\caption{Largest-angle bisection. The largest angle \(\alpha\) is bisected, producing the two child triangles \(A_L\) and \(A_R\).}
\label{fig:lab-operation}
\end{figure}

We now describe LAB in the angle coordinates introduced above. Let the parent triangle have ordered angles \((\alpha,\beta,\gamma)\), so that \(\alpha\) is the angle selected for bisection. As shown in Figure~\ref{fig:lab-operation}, the bisector divides \(\alpha\) into two angles of size \(\alpha/2\). The third angle of each child is determined by the angle sum. Before reordering, the two children therefore have angle triples
\begin{equation}\label{eq:rawchildren}
A_L(\alpha,\beta,\gamma)
=
\left(\frac{\alpha}{2},\beta,\gamma+\frac{\alpha}{2}\right),
\qquad
A_R(\alpha,\beta,\gamma)
=
\left(\frac{\alpha}{2},\gamma,\beta+\frac{\alpha}{2}\right).
\end{equation}
We call these triples \emph{raw} because their entries are inherited from the geometry of the two children and need not be arranged in decreasing order.

Our shape chamber always records angles in decreasing order. After bisection, however, the largest child angle need not occupy the same coordinate as the largest parent angle. We must therefore sort each raw child triple before applying the projective chart. If \(\sortop\) denotes decreasing rearrangement, the normalized LAB maps are
\begin{equation}\label{eq:LRdefinition}
L=\chi\circ\sortop\circ A_L\circ\chi^{-1},
\qquad
R=\chi\circ\sortop\circ A_R\circ\chi^{-1}.
\end{equation}
The labels \(L\) and \(R\) refer to the two children shown in Figure~\ref{fig:lab-operation}. Because the raw maps preserve positivity and the angle sum, sorting returns each child to \(\D\). Consequently, \(L\) and \(R\) are self-maps of \(\T\).

\subsection{Piecewise formulas}

In the chart \(\T\), the raw triples are projectively equivalent to
\[
\left(\frac12,x,y+\frac12\right)
\qquad\text{and}\qquad
\left(\frac12,y,x+\frac12\right).
\]
The projective expressions themselves are simple. Their piecewise character comes entirely from determining the decreasing order of the three entries. The boundaries between the cases are therefore sorting walls, where two raw child angles coincide.

\begin{proposition}[Piecewise projective LAB maps]\label{prop:formulas}
The left map is
\begin{equation}\label{eq:Lformula}
L(x,y)=
\begin{cases}
\displaystyle
\left(\frac{x}{y+1/2},\frac{1}{2(y+1/2)}\right),
& x\geq\frac12,\ y\geq x-\frac12,\\[10pt]
\displaystyle
\left(\frac{y+1/2}{x},\frac{1}{2x}\right),
& x\geq\frac12,\ y\leq x-\frac12,\\[10pt]
\displaystyle
\left(\frac{1}{2(y+1/2)},\frac{x}{y+1/2}\right),
& x\leq\frac12.
\end{cases}
\end{equation}
The right map is
\begin{equation}\label{eq:Rformula}
R(x,y)=
\begin{cases}
\displaystyle
\left(\frac{y}{x+1/2},\frac{1}{2(x+1/2)}\right),
& y\geq\frac12,\\[10pt]
\displaystyle
\left(\frac{1}{2(x+1/2)},\frac{y}{x+1/2}\right),
& y\leq\frac12.
\end{cases}
\end{equation}
In particular, \(L\) and \(R\) are continuous piecewise projective self-maps of \(\T\).
\end{proposition}

\begin{proof}
Consider first the raw left-child coordinates
\[
\left(\frac12,x,y+\frac12\right).
\]
Since \(y+1/2\geq1/2\), the largest entry is either \(x\) or \(y+1/2\). If \(x\geq1/2\) and \(y\geq x-1/2\), then
\[
y+\frac12\geq x\geq\frac12.
\]
After sorting and dividing the two smaller entries by the largest one, we obtain the first branch of \eqref{eq:Lformula}. If \(x\geq1/2\) and \(y\leq x-1/2\), then
\[
x\geq y+\frac12\geq\frac12,
\]
which gives the second branch. Finally, if \(x\leq1/2\), then the defining inequalities \(x\geq y\) imply
\[
y+\frac12\geq\frac12\geq x,
\]
and sorting gives the third branch.

For the right child, the raw coordinates are
\[
\left(\frac12,y,x+\frac12\right).
\]
Here \(x+1/2\) is always the largest entry, and the remaining two entries exchange order at \(y=1/2\). This gives the two branches of \eqref{eq:Rformula}. On each sorting wall, the adjacent expressions agree because the exchanged coordinates are equal. Hence both maps are continuous.
\end{proof}

\section{Dynamics along the branches of the refinement tree}\label{sec:dynamics}

The purpose of this section is to classify the asymptotic behavior of an infinite LAB branch. Here we show that the two child maps are individually non-expansive with respect to a projective distance from the right-isosceles shape, but neither map is strictly contractive: each fixes an entire segment of shapes. These fixed segments describe the exceptional behavior that can persist when the same child is selected repeatedly.

The situation changes when the selected child changes. We show that every mixed two-step block decreases the projective distance strictly, except at the right-isosceles shape. The proof therefore has three stages. We first record the compactness inherited from deterministic angle regularity and identify the fixed sets of the pure maps. We then introduce a projective Lyapunov function and determine precisely where one-step decrease can fail. Finally, we combine these facts to describe pure branches and prove that every branch with infinitely many changes of child converges to the unique intersection of the two pure fixed sets.

\subsection{Compactness and pure fixed sets}

We use the known deterministic regularity of LAB through the following compactness consequence.

\begin{lemma}[Compactness of LAB orbits]\label{lem:orbitcompactness}
For every initial shape $z_0\in\T$, there exists $\eta>0$ such that every descendant of $z_0$ belongs to
\[
K_\eta=\{(x,y)\in\T:y\geq\eta\}.
\]
In particular, every infinite LAB branch issued from $z_0$ has compact closure in $\T$.
\end{lemma}

\begin{proof}
Let $\alpha_0\geq\beta_0\geq\gamma_0>0$ be the ordered angles of the initial triangle. Ismailescu et al.~\cite{ismailescu2019largestanglebisectionprocedure} prove that every angle generated by LAB is bounded below by
\[
c_0=\min\left\{\gamma_0,\frac{\alpha_0}{2}\right\}>0\,,
\]
see also \cite{michaudKorotov2026}. If $(\alpha,\beta,\gamma)$ is any descendant angle triple, then $\gamma\geq c_0$ and $\alpha\leq\pi$. Consequently,
\[
y=\frac{\gamma}{\alpha}\geq\frac{c_0}{\pi}.
\]
Thus the orbit is contained in $K_\eta$ with $\eta=c_0/\pi$. Since
\[
K_\eta=\{(x,y):1\geq x\geq y\geq\eta\}
\]
is compact, the final assertion follows.
\end{proof}

The exceptional pure dynamics are already visible in the fixed sets of the two child maps.

\begin{proposition}[Pure fixed sets]\label{prop:fixedsets}
The fixed sets of the two LAB maps are
\begin{equation}\label{eq:CLCR}
C_L=\Fix(L)=\{(x,1/2):1/2\leq x\leq1\},
\qquad
C_R=\Fix(R)=\{(1/2,y):0<y\leq1/2\}.
\end{equation}
Their intersection is the single point
\begin{equation}\label{eq:p}
p=\left(\frac12,\frac12\right),
\end{equation}
which represents the right-isosceles angle triple $(\pi/2,\pi/4,\pi/4)$.
\end{proposition}

\begin{proof}
For $1/2\leq x\leq1$, the first branch of \eqref{eq:Lformula} gives
\[
L(x,1/2)=(x,1/2),
\]
so the displayed horizontal segment is contained in $\Fix(L)$. Conversely, solving $L(x,y)=(x,y)$ in each of the three regions of \eqref{eq:Lformula} yields $y=1/2$ and $1/2\leq x\leq1$. Hence this segment is the whole fixed set of $L$.

Similarly, the formulas in \eqref{eq:Rformula} show that $R(x,y)=(x,y)$ exactly when $x=1/2$ and $0<y\leq1/2$. The two fixed segments meet only at $p$.
\end{proof}

Thus neither child map has a unique attractor under pure iteration. Each admits a continuum of stationary shapes, and the two continua meet only at the right-isosceles point. This already suggests that convergence to $p$ cannot be produced by either map alone; it must arise from changing between them.

\subsection{A Hilbert-projective Lyapunov function}

Because the shape domain $\T$ is obtained by projective normalization, it is natural to measure distances using Hilbert's projective metric. For positive triples $u,v$, this metric is defined by
\[
d_H(u,v)
=
\log\frac{\max_i(u_i/v_i)}{\min_i(u_i/v_i)}.
\]
For background on Hilbert's projective metric and non-expansive maps, see \cite{Nussbaum1988}.

We use the same notation for the induced metric on $\T$. If $z=(x,y)$ and $w=(x',y')$, we set
\[
d_H(z,w)
=
d_H\bigl((1,x,y),(1,x',y')\bigr).
\]
This definition is independent of the chosen positive homogeneous representatives because Hilbert's metric is invariant under positive rescaling. For a nonempty subset $C\subset\T$, we write
\[
d_H(z,C)=\inf_{w\in C}d_H(z,w).
\]

Recall that
\[
p=\left(\frac12,\frac12\right)
\]
represents the right-isosceles shape. We define the Lyapunov function
\[
D(z)=d_H(z,p).
\]
For $z=(x,y)$, the homogeneous representatives
\[
(1,x,y)
\qquad\text{and}\qquad
\left(1,\frac12,\frac12\right)
\]
give the coordinate ratios $1$, $2x$, and $2y$. Since $x\geq y$, it follows that
\begin{equation}\label{eq:D}
D(x,y)=\log H(x,y),
\qquad
H(x,y)=\frac{\max\{x,1/2\}}{\min\{y,1/2\}}.
\end{equation}
In particular, $D\geq0$, and $D(z)=0$ if and only if $z=p$.

On every compact subset of $\T$ bounded away from the degenerate boundary $y=0$, the Hilbert metric induces the usual coordinate topology. Lemma~\ref{lem:orbitcompactness} therefore allows us to use compactness and convergence in either description along a LAB orbit.

The next lemma shows that $D$ is a Lyapunov function for either child map. The proof uses the homogeneous angle triples directly and does not require the piecewise formulas of Proposition~\ref{prop:formulas}.

\begin{lemma}[One-step non-expansiveness]\label{lem:nonexpansive}
For every $z\in\T$,
\[
D(Lz)\leq D(z),
\qquad
D(Rz)\leq D(z).
\]
\end{lemma}

\begin{proof}
Let $u=(\alpha,\beta,\gamma)\in\D$ and set
\[
m=\min\left\{\frac\alpha2,\beta,\gamma\right\},
\qquad
M=\max\left\{\frac\alpha2,\beta,\gamma\right\}.
\]
The point $p$ has homogeneous angle proportions $(2,1,1)$. Hence the coordinate ratios entering $d_H(u,p)$ are
\[
\frac{\alpha}{2},\qquad \beta,\qquad \gamma,
\]
and therefore
\[
D(u)=d_H(u,p)=\log\frac{M}{m}.
\]

Consider the raw left child
\[
\left(\frac\alpha2,\beta,\gamma+\frac\alpha2\right),
\]
and let $(\alpha',\beta',\gamma')$ be its decreasing rearrangement. We verify that each of the three quantities
\[
\frac{\alpha'}2,\qquad \beta',\qquad \gamma'
\]
belongs to $[m,M]$.

First, every raw child angle is at least $m$. The only raw entry that can exceed $M$ is $\gamma+\alpha/2$, and if it does, it becomes the largest sorted angle. Hence the two smaller sorted angles $\beta'$ and $\gamma'$ lie in $[m,M]$. Second,
\[
2m\leq\gamma+\frac\alpha2\leq2M,
\]
while the other two raw entries already lie in $[m,M]$. It follows that the largest sorted angle satisfies
\[
m\leq\frac{\alpha'}2\leq M.
\]
Thus the spread of the three ratios defining the projective distance of the child is no greater than $M/m$, and therefore
\[
D(Lu)\leq\log\frac{M}{m}=D(u).
\]
The right child is treated identically, with $\beta+\alpha/2$ in place of $\gamma+\alpha/2$.
\end{proof}

Non-expansiveness alone does not imply convergence: the value of $D$ may remain constant along nontrivial orbits. We must therefore determine exactly where a single LAB step fails to decrease it.

\begin{lemma}[Equality sets]\label{lem:equalitysets}
Equality in Lemma~\ref{lem:nonexpansive} holds for $L$ exactly on
\begin{equation}\label{eq:EL}
E_L
=
\{(x,y)\in\T:y\geq1/2\}
\cup
\{(x,y)\in\T:x=y\leq1/2\},
\end{equation}
and for $R$ exactly on
\begin{equation}\label{eq:ER}
E_R
=
\{(x,y)\in\T:x\leq1/2\}
\cup
\{(x,y)\in\T:x=y\geq1/2\}.
\end{equation}
\end{lemma}

\begin{proof}
We work with $H=e^D$ and use the formulas of Proposition~\ref{prop:formulas}.

\smallskip
\noindent\emph{Equality for $L$.}
There are three cases.

If $y\geq1/2$, then $H(x,y)=2x$. The first formula in \eqref{eq:Lformula}, written $L(x,y)=(X,Y)$, satisfies $X\geq1/2$, $Y\leq1/2$, and
\[
\frac{X}{Y}=2x.
\]
Hence $H(L(x,y))=H(x,y)$ throughout this region.

Suppose next that $x\geq1/2>y$. Then $H(x,y)=x/y$. In the first $L$-region,
\[
H(L(x,y))=\frac{2x}{y+1/2}<\frac{x}{y},
\]
because $y<1/2$. In the second $L$-region,
\[
H(L(x,y))=\frac{2y+1}{x}.
\]
The region condition $x\geq y+1/2$, together with $0<y<1/2$, gives
\[
x^2\geq\left(y+\frac12\right)^2>2y^2+y.
\]
Equivalently,
\[
\frac{2y+1}{x}<\frac{x}{y}.
\]
Thus the decrease is strict throughout this case.

Finally, suppose that $x\leq1/2$. Then $H(x,y)=1/(2y)$. The third formula in \eqref{eq:Lformula} gives
\[
H(L(x,y))\leq\frac{1}{2x}\leq\frac{1}{2y}.
\]
Equality holds precisely when $x=y$. Hence the equality set for $L$ is exactly $E_L$.

\smallskip
\noindent\emph{Equality for $R$.}
If $x\leq1/2$, then $H(x,y)=1/(2y)$. Using the second formula in \eqref{eq:Rformula}, one obtains $H(R(x,y))=H(x,y)$. If $x\geq1/2$ and $x>y$, direct substitution in the appropriate branch of~\eqref{eq:Rformula} gives strict decrease. On the remaining diagonal $x=y\geq1/2$, equality holds. This yields exactly the set $E_R$ in \eqref{eq:ER}.
\end{proof}

\subsection{Asymptotic behavior of pure branches}

The equality sets $E_L$ and $E_R$ are larger than the fixed segments $C_L$ and $C_R$. A point outside a fixed segment can therefore preserve $D$ for one step. The following lemma shows that such equality cannot persist under a second application of the same map. Persistent failure of contraction under pure iteration is possible only on the corresponding fixed set.

\begin{lemma}[Invariant equality sets]\label{lem:invariantequality}
One has
\[
\{z\in E_L:Lz\in E_L\}=C_L,
\qquad
\{z\in E_R:Rz\in E_R\}=C_R.
\]
\end{lemma}

\begin{proof}
We give the argument for $L$. Suppose first that $z=(x,y)\in E_L$ with $y\geq1/2$. Writing $Lz=(X,Y)$, the first formula of \eqref{eq:Lformula} gives
\[
Y=\frac{1}{2(y+1/2)}\leq\frac12,
\]
with equality only when $y=1/2$. If $Lz\in E_L$, then either $Y\geq1/2$, which forces $y=1/2$, or $X=Y\leq1/2$, which forces $x=1/2$ and hence $x=y=1/2$. In either case, $z\in C_L$.

Now suppose that $z=(x,x)$ with $x\leq1/2$. Then
\[
L(x,x)
=
\left(\frac{1}{2x+1},\frac{x}{x+1/2}\right).
\]
Requiring this image to belong to $E_L$ forces $x=1/2$, again giving a point of $C_L$. Conversely, $L$ fixes $C_L$ pointwise. Hence the first identity follows. The second is obtained analogously from~\eqref{eq:Rformula}.
\end{proof}

We can now describe the limiting behavior of a branch that repeatedly chooses the same child. The conclusion is attraction to a fixed segment.

\begin{proposition}[Attraction to the pure fixed sets]\label{prop:purestability}
For every $z\in\T$,
\[
d_H(L^n z,C_L)\longrightarrow0,
\qquad
d_H(R^n z,C_R)\longrightarrow0.
\]
\end{proposition}

\begin{proof}
We prove the assertion for $L$. By Lemma~\ref{lem:orbitcompactness}, the orbit $z_n=L^nz$ has compact closure in $\T$. By Lemma~\ref{lem:nonexpansive}, the sequence $D(z_n)$ is non-increasing and therefore converges to some $\ell\geq0$.

Let $z_*$ be an accumulation point, and choose a subsequence $z_{n_k}\to z_*$. Continuity of $D$ and $L$ gives
\[
D(z_*)=\ell,
\qquad
D(Lz_*)=\lim_{k\to\infty}D(z_{n_k+1})=\ell.
\]
Hence $z_*\in E_L$. Moreover, $Lz_*$ is itself an accumulation point of the orbit, so the same argument gives $Lz_*\in E_L$. Lemma~\ref{lem:invariantequality} therefore implies $z_*\in C_L$.

Suppose that $d_H(L^nz,C_L)$ does not tend to zero. Then there exist $\varepsilon>0$ and a subsequence $z_{n_k}$ such that
\[
d_H(z_{n_k},C_L)\geq\varepsilon
\]
for every $k$. Compactness provides a further subsequence converging to an accumulation point $z_*$. As shown above, $z_*\in C_L$. Since the Hilbert metric induces the coordinate topology on the compact set containing the orbit,
\[
d_H(z_{n_k},C_L)
\leq d_H(z_{n_k},z_*)\longrightarrow0,
\]
a contradiction. Therefore
\[
d_H(L^nz,C_L)\longrightarrow0.
\]
The proof for $R$ is identical.
\end{proof}

\subsection{Mixed branches and global classification}

The preceding results explain why one child map alone cannot force convergence to $p$. We now show that changing child removes this obstruction. Although neither $L$ nor $R$ strictly decreases $D$ everywhere, each mixed two-step block does so away from the common fixed point.

\begin{theorem}[Strict Lyapunov decrease for mixed blocks]\label{thm:mixedblock}
For every $z\in\T\setminus\{p\}$,
\[
D(RLz)<D(z),
\qquad
D(LRz)<D(z).
\]
\end{theorem}

\begin{proof}
By Lemma~\ref{lem:nonexpansive},
\[
D(RLz)\leq D(Lz)\leq D(z).
\]
If equality held throughout, Lemma~\ref{lem:equalitysets} would imply
\[
z\in E_L
\qquad\text{and}\qquad
Lz\in E_R.
\]
We show that these two conditions force $z=p$.

Suppose first that $y\geq1/2$. Then
\[
L(x,y)
=
\left(\frac{x}{y+1/2},\frac{1}{2(y+1/2)}\right)
=(X,Y),
\qquad
Y\leq\frac12.
\]
For $(X,Y)$ to belong to $E_R$, either $X\leq1/2$ or $X=Y=1/2$. In the first case,
\[
2x\leq y+\frac12.
\]
Since $x\geq y\geq1/2$, this inequality forces $x=y=1/2$. The second case gives the same conclusion.

The remaining component of $E_L$ is the diagonal $x=y\leq1/2$. There,
\[
L(x,x)
=
\left(\frac{1}{2x+1},\frac{x}{x+1/2}\right).
\]
Requiring this point to lie in $E_R$ again forces $x=1/2$. Thus equality in the chain can occur only at $p$, proving
\[
D(RLz)<D(z)
\]
for every $z\neq p$.

The proof for $LR$ is symmetric. Equality would require $z\in E_R$ and $Rz\in E_L$, and the two components of $E_R$, together with \eqref{eq:Rformula}, again force $z=p$.
\end{proof}

Strict decrease at each change of child is the essential local mechanism. To obtain a global convergence theorem, we must also rule out the possibility that the amount of decrease tends to zero while the orbit remains away from $p$. Compactness provides exactly this uniformity: on any compact set that excludes $p$, both mixed blocks decrease $D$ by a common positive amount.

\begin{theorem}[Classification of infinite branches]\label{thm:classification}
Let $z_0\in\T$, let $\omega_n\in\{L,R\}$, and define
\[
z_{n+1}=\omega_{n+1}(z_n).
\]
If the infinite word $\omega_1\omega_2\cdots$ has infinitely many changes of symbol, then
\[
z_n\longrightarrow p.
\]
Consequently, every infinite LAB branch is of exactly one of the following two types:
\begin{enumerate}
\item the itinerary is not eventually constant, in which case the shapes converge to $p$;
\item the itinerary is eventually constant, in which case the branch approaches $C_L$ if its eventual symbol is $L$, and approaches $C_R$ if its eventual symbol is $R$.
\end{enumerate}
\end{theorem}

\begin{proof}
By Lemma~\ref{lem:nonexpansive}, the sequence $D(z_n)$ is non-increasing, and hence converges to some $\ell\geq0$. Suppose, for contradiction, that $\ell>0$. Choose $\varepsilon$ with
\[
0<\varepsilon<\ell.
\]
By Lemma~\ref{lem:orbitcompactness}, the orbit is contained in $K_\eta$ for some $\eta>0$. The set
\[
K_{\eta,\varepsilon}
=
\{z\in K_\eta:D(z)\geq\varepsilon\}
\]
is compact and does not contain $p$.

By Theorem~\ref{thm:mixedblock}, the continuous functions
\[
z\longmapsto D(z)-D(RLz),
\qquad
z\longmapsto D(z)-D(LRz)
\]
are strictly positive on $K_{\eta,\varepsilon}$. They therefore have a common positive lower bound $\delta_\varepsilon$ on that set.

Infinitely many changes of symbol yield infinitely many pairwise disjoint blocks $LR$ or $RL$. At the beginning of each such block, the orbit lies in $K_{\eta,\varepsilon}$ because
\[
D(z_n)\geq\ell>\varepsilon.
\]
Every mixed block therefore lowers $D$ by at least $\delta_\varepsilon$, while all intervening steps are non-increasing. After sufficiently many mixed blocks, this would force $D$ to become negative, contradicting $D\geq0$. Hence $\ell=0$.

Since $D$ vanishes only at $p$, it remains to convert $D(z_n)\to0$ into convergence of the shapes. The orbit lies in the compact set $K_\eta$, and every accumulation point $z_*$ satisfies $D(z_*)=0$ by continuity. Thus every accumulation point equals $p$, and consequently $z_n\to p$.

Finally, an infinite binary word has only finitely many changes if and only if it is eventually constant. The second alternative then follows from Proposition~\ref{prop:purestability}.
\end{proof}

\section{Statistical regularization}\label{sec:statistics}

The branch classification obtained in Section~\ref{sec:dynamics} has two closely related global consequences. First, if the two children are selected randomly with positive probability, then a typical branch changes symbol infinitely often and therefore converges to the right-isosceles shape. Second, under uniform child selection, a random word of length $n$ samples one of the $2^n$ level-$n$ descendants uniformly. Almost-sure convergence along random branches can thus be translated into concentration of the complete refinement tree near the right-isosceles shape.

\subsection{Random branches}

Let $(\omega_n)_{n\geq1}$ be independent random variables taking values in $\{L,R\}$. Starting from $z_0\in\T$, define the corresponding random branch by
\[
z_n(\omega)
=
\omega_n\circ\cdots\circ\omega_1(z_0).
\]

\begin{theorem}[Almost-sure regularization]\label{thm:as}
Let $0<q<1$, and suppose that
\[
\mathbb P(\omega_n=L)=q,
\qquad
\mathbb P(\omega_n=R)=1-q.
\]
Then, for every initial shape $z_0\in\T$,
\[
z_n(\omega)\longrightarrow p
\]
for almost every itinerary $\omega$.
\end{theorem}

\begin{proof}
For $N\geq1$, let
\[
A_N^L
=
\{\omega_n=L\text{ for every }n\geq N\}.
\]
By independence,
\[
\mathbb P(A_N^L)
=
\lim_{k\to\infty}q^k
=
0.
\]
Similarly, the event
\[
A_N^R
=
\{\omega_n=R\text{ for every }n\geq N\}
\]
has probability zero. The event that the itinerary is eventually constant is
\[
\bigcup_{N\geq1}\bigl(A_N^L\cup A_N^R\bigr),
\]
and therefore also has probability zero. Almost every itinerary has infinitely many changes of symbol, so Theorem~\ref{thm:classification} applies and gives $z_n(\omega)\to p$.
\end{proof}

The restriction $0<q<1$ ensures that neither child is excluded. The limiting shape does not otherwise depend on the selection bias: changing $q$ changes the relative frequency of the two symbols, but almost every itinerary still changes symbol infinitely often.

\subsection{Complete refinement levels}

We now specialize to uniform child selection, $q=1/2$. Every word of length $n$ then has probability $2^{-n}$, so the law of the random shape $z_n$ is precisely the empirical distribution of the shapes of all $2^n$ descendants at refinement level $n$.

For a word $w=\omega_n\cdots\omega_1\in\{L,R\}^n$, write
\[
w(z_0)
=
\omega_n\circ\cdots\circ\omega_1(z_0).
\]
The empirical shape measure of the complete level-$n$ refinement tree is
\begin{equation}\label{eq:mun}
\mu_n
=
2^{-n}\sum_{w\in\{L,R\}^n}\delta_{w(z_0)}.
\end{equation}
Each descendant is therefore counted with equal weight, including multiplicities when distinct words produce the same normalized shape.

\begin{corollary}[Statistical regularization of the refinement tree]\label{cor:tree}
For every initial shape $z_0\in\T$,
\[
\mu_n\Rightarrow\delta_p.
\]
Equivalently, for every neighborhood $U$ of $p$,
\[
\frac{\#\{w\in\{L,R\}^n:w(z_0)\in U\}}{2^n}
\longrightarrow1.
\]
\end{corollary}

\begin{proof}
Under independent uniform choices of $L$ and $R$, the law of $z_n$ is exactly $\mu_n$. By Theorem~\ref{thm:as},
\[
z_n\longrightarrow p
\qquad\text{almost surely}.
\]
Let $f:\T\to\mathbb R$ be bounded and continuous. Then
\[
f(z_n)\longrightarrow f(p)
\qquad\text{almost surely},
\]
and bounded convergence gives
\[
\int_{\T}f\,d\mu_n
=
\mathbb E[f(z_n)]
\longrightarrow
f(p)
=
\int_{\T}f\,d\delta_p.
\]
Hence $\mu_n\Rightarrow\delta_p$.

Now let $U$ be a neighborhood of $p$. Since the limiting measure is concentrated at $p$, weak convergence implies $\mu_n(U)\to1$. By the definition of $\mu_n$, this is precisely the stated proportion of level-$n$ descendants.
\end{proof}

Corollary~\ref{cor:tree} is the population-level consequence of the branch classification. It does not assert that every branch converges to $p$: eventually pure branches retain the exceptional behavior described in Proposition~\ref{prop:purestability}. Rather, it shows that this exceptional behavior is asymptotically negligible when all descendants at a fixed refinement level are assigned equal weight. In this precise sense, the right-isosceles shape statistically dominates the complete refinement tree.

\section{Conclusion}\label{sec:conclusion}

Largest-angle bisection admits a natural description as a two-generator piecewise projective system on triangle shape space. This viewpoint reveals a sharp distinction between pure and mixed branches. Repeated selection of a single child drives the dynamics toward the corresponding one-dimensional fixed segment, whereas every itinerary with infinitely many changes converges to the right-isosceles shape. Thus LAB does not drive every infinite branch to a single shape, but it does possess a unique statistically dominant shape.

The branchwise classification has a global consequence. Since eventually pure itineraries form a null set for every Bernoulli process assigning positive probability to both children, a random branch converges to the right-isosceles shape almost surely. Under uniform selection, this becomes a statement about the complete refinement tree: the proportion of level-$n$ descendants outside any prescribed neighborhood of the right-isosceles shape tends to zero. The result therefore goes beyond standard deterministic shape regularity. A positive lower bound on the angles prevents degeneration, whereas the present analysis identifies the geometry that dominates the refinement tree. This dynamics is rather different than what happens for LEB. Kalmanovich et al.~\cite{kalmanovich2026stability} showed that LEB branches eventually cycle with a cycle length of length at most four and that these cycles belong to a compact subspace of the triangle space. Since every initial triangles yields a different limiting cycle, a dense subspace of the space of triangle shapes can be reached. This is rather different that almost sure convergence to a single shape.

From the perspective of finite element mesh refinement, the result should be interpreted as structural rather than immediately algorithmic. The present analysis concerns the unconstrained binary refinement tree generated from a single triangle and assigns equal weight to every descendant at a fixed depth. It does not yet incorporate conforming completion, interactions between neighboring elements, adaptive marking, or weights based on area or error contribution. Nevertheless, LAB provides a simple setting in which geometric refinement, projective dynamics, and probability interact transparently, and in which a distributional statement can be obtained in addition to worst-case angle bounds.

Several questions remain open. On the quantitative side, it would be useful to obtain explicit convergence estimates in terms of the number and spacing of changes between the two children, together with finite-level bounds for the proportion of descendants outside a prescribed neighborhood of the right-isosceles shape. On the algorithmic side, one may ask whether comparable statistical regularization persists under conforming completion, adaptive refinement, or area- and error-weighted sampling. On the geometric side, extensions to largest-angle $n$-section and to higher-dimensional refinement rules based on solid or dihedral angles merit investigation. These questions connect the present geometric analysis with the broader study of robust mesh generation and adaptive refinement for finite element computation.

\bibliography{references}

\end{document}